\documentclass[a4paper,12pt]{amsart}
\usepackage{amsfonts}
\usepackage{amssymb}
\usepackage{ifthen}
\usepackage{graphicx}
\nonstopmode \numberwithin{equation}{section}
\usepackage[usenames]{color}
\newtheorem{thm}{Theorem}
\newtheorem{lem}{Lemma}
\newtheorem{cor}{Corollary}
\newtheorem{prop}{Proposition}

\newtheorem{cl}{Claim}
\newtheorem{ca}{Case}
\newtheorem{sca}{Subcase}
\newtheorem{scl}{Subclaim}
\newtheorem{conj}[equation]{Conjecture}
\newtheorem{examp}{Example}
\theoremstyle{definition}
\newtheorem{defn}{Definition}

\newtheorem{op}[equation]{Open Problem}
\newtheorem{ques}[equation]{Question}
\newtheorem{rem}{Remark}
\newtheorem{exam}[equation]{Example}

\newcounter {own}
\def\theown {\thesection       .\arabic{own}}

\newenvironment{pf}[1][]{%
 \vskip 3mm
 \noindent
 \ifthenelse{\equal{#1}{}}%
  {{\slshape Proof. }}%
  {{\slshape #1.} }%
 }%
{\qed\bigskip}

\newcounter{alphabet}
\newcounter{tmp}
\newenvironment{Thm}[1][]{\refstepcounter{alphabet}%
\bigskip%
\noindent%
{\bf Theorem \Alph{alphabet}}%
\ifthenelse{\equal{#1}{}}{}{ (#1)}%
{\bf .} \itshape}{\vskip 8pt}

\makeatletter
\newcommand{\Ref}[1]{\@ifundefined{r@#1}{}{\setcounter{tmp}{\ref{#1}}\Alph{tmp}}}
\makeatother

\newcommand{\IC}{{\mathbb C}}

\def\be{\begin{equation}}
\def\ee{\end{equation}}

\newcommand{\bee}{\begin{enumerate}}
\newcommand{\eee}{\end{enumerate}}

\newcommand{\blem}{\begin{lem}}
\newcommand{\elem}{\end{lem}}
\newcommand{\bthm}{\begin{thm}}
\newcommand{\ethm}{\end{thm}}
\newcommand{\bcor}{\begin{cor}}
\newcommand{\ecor}{\end{cor}}
\newcommand{\beg}{\begin{exam}}
\newcommand{\eeg}{\end{exam}}
\newcommand{\begs}{\begin{examples}}
\newcommand{\eegs}{\end{examples}}
\newcommand{\bdefe}{\begin{defn}}
\newcommand{\edefe}{\end{defn}}
\newcommand{\bprob}{\begin{prob}}
\newcommand{\eprob}{\end{prob}}
\newcommand{\bques}{\begin{ques}}
\newcommand{\eques}{\end{ques}}
\newcommand{\bei}{\begin{itemize}}
\newcommand{\eei}{\end{itemize}}
\newcommand{\bcon}{\begin{conj}}
\newcommand{\econ}{\end{conj}}
\newcommand{\bop}{\begin{op}}
\newcommand{\eop}{\end{op}}
\newcommand{\bca}{\begin{ca}}
\newcommand{\eca}{\end{ca}}
\newcommand{\bsca}{\begin{sca}}
\newcommand{\esca}{\end{sca}}

\newcommand{\bcl}{\begin{cl}}
\newcommand{\ecl}{\end{cl}}

\newcommand{\bscl}{\begin{scl}}
\newcommand{\escl}{\end{scl}}

\newcommand{\bcons}{\begin{conjs}}
\newcommand{\econs}{\end{conjs}}
\newcommand{\bprop}{\begin{propo}}
\newcommand{\eprop}{\end{propo}}
\newcommand{\br}{\begin{rem}}
\newcommand{\er}{\end{rem}}
\newcommand{\brs}{\begin{rems}}
\newcommand{\ers}{\end{rems}}
\newcommand{\bo}{\begin{obser}}
\newcommand{\eo}{\end{obser}}
\newcommand{\bos}{\begin{obsers}}
\newcommand{\eos}{\end{obsers}}
\newcommand{\bpf}{\begin{pf}}
\newcommand{\epf}{\end{pf}}
\newcommand{\ba}{\begin{array}}
\newcommand{\ea}{\end{array}}
\newcommand{\beq}{\begin{eqnarray}}
\newcommand{\beqq}{\begin{eqnarray*}}
\newcommand{\eeq}{\end{eqnarray}}
\newcommand{\eeqq}{\end{eqnarray*}}

\newcounter{minutes}
\divide\time by 60
\newcounter{hours}
\multiply\time by 60 \addtocounter{minutes}{-\time}

\begin{document}
\bibliographystyle{amsplain}
\title [] {Starlikeness and convexity of polyharmonic mappings}

\def\thefootnote{}
\footnotetext{ \texttt{\tiny File:~\jobname .tex,
          printed: \number\day-\number\month-\number\year,
          \thehours.\ifnum\theminutes<10{0}\fi\theminutes}
} \makeatletter\def\thefootnote{\@arabic\c@footnote}\makeatother

\author{J. Chen}
\address{J. Chen, Department of Mathematics,
Hunan Normal University, Changsha, Hunan 410081, People's Republic
of China.} \email{jiaolongchen@sina.com}

\author{ A. Rasila }
\address{A. Rasila, Department of Mathematics,
Hunan Normal University, Changsha, Hunan 410081, People's Republic
of China, and
Department of Mathematics and Systems Analysis, Aalto University, P.O. Box 11100, FI-00076 Aalto,
 Finland.} \email{antti.rasila@iki.fi}

\author{X. Wang${}^{~\mathbf{*}}$}
\address{X. Wang, Department of Mathematics,
Hunan Normal University, Changsha, Hunan 410081, People's Republic
of China.} \email{xtwang@hunnu.edu.cn}


\begin{abstract}
In this paper, we first find an estimate for the range of polyharmonic mappings in the class $HC_{p}^{0}$. Then, we obtain two characterizations in terms of the convolution for polyharmonic mappings to be starlike of order $\alpha$, and convex of order $\beta$, respectively. Finally, we study the radii of starlikeness and convexity for polyharmonic mappings, under certain coefficient conditions.
\end{abstract}


\keywords{polyharmonic mapping, starlike,
 convex\\
${}^{\mathbf{*}}$ Corresponding author}

\maketitle

\section{Introduction }\label{csw-sec1}

Let $F=u+iv$ be a $2p$ times continuously differentiable complex-valued mapping, where $p\geq1$, defined in a domain $D \subset \mathbb{C}$. The mapping $F$ is called {\it polyharmonic} (or {\it $p$-harmonic}), if it satisfies the polyharmonic equation $\Delta^{p}F =\Delta(\Delta^{p-1}F)= 0$, where $\Delta^{1}:=\Delta$ is the usual complex Laplacian operator
$$
\Delta=4\frac{\partial^{2}}{\partial z\partial \overline{z}}:=
\frac{\partial^{2}}{\partial x^{2}}+\frac{\partial^{2}}{\partial y^{2}}.
$$
For a simply connected domain $D$, it is well known (see \cite{CPW2, qiwa}) that  a mapping $F$ is polyharmonic if and only if $F$  has the representation
$$F(z)=\sum_{k=1}^{p}|z|^{2(k-1)}G_{k}(z),$$ where $G_{k}$ are
complex-valued harmonic mappings in $D$ for all $k\in \{1,\cdots,p\}$. The mappings $G_{k}$ can be presentented (see \cite{cl, du}) in the form
$$
G_{k} = h_{k} + \overline{g_{k}},
$$
where all $h_{k}$ and  $g_{k}$ are analytic in $D$ for all $k\in \{1,\cdots,p\}$. Clearly, if $p=1$, then we have the usual class of harmonic mappings and, for $p=2$, we obtain the class of biharmonic mappings, as special cases.

The biharmonic equation is related to numerous modeling problems in science and engineering. For example, it arises from certain problems in solid mechanics, and also from the theory of steady Stokes flow (i.e., speed $\approx 0$) of viscous fluids, where it is the equation satisfed by the stream function (see e.g. \cite{ha, kh, la}). In the geometric function theory, the class biharmonic mappings can be understood as a natural generalization of the harmonic mappings, but it has only recently been studied from this point of view (see \cite{abab, A, ababkh, CPW0, CPW4,CW}). The reader is referred to \cite{CRW, CPW2, qiwa} for the properties of polyharmonic mappings, and \cite{cl, du} for basic results on harmonic mappings.

For $r>0$, write $\mathbb{U}_{r} =\{z:|z|<r\}$ $(r >0$), and let $\mathbb{U}:=\mathbb{U}_{1}$, i.e., the unit disk. Let $S_{H}$ denote the set of all univalent harmonic mappings $f$ in $\mathbb{U}$, where
\be\label{eq1.2}f(z)=h(z)+\overline{g(z)}=z+\sum_{j=2}^{\infty}a_j
z^j+\sum_{j=1}^{\infty}\overline{b_jz^j}\ee with $|b_{1}|<1$. We denote by $S_{H}^{0}$ the set of all mappings in
$S_{H}$ with $b_{1}=0$. Let $S^{\ast}_{H}$ and $S^{\ast,0}_{H}$
denote the respective subclasses of $S_{H}$ and $S^{0}_{H}$, where the images of
$f(\mathbb{U})$ are starlike. Let $C_{H}$ and $C^{0}_{H}$
denote the respective subclasses of $S_{H}$ and $S^{0}_{H}$, where the images of
$f(\mathbb{U})$ are convex.

In \cite{av}, Avci and Z{\l}otkiewicz introduced the class $HS$ of univalent harmonic mappings $F$ with the series expansion \eqref{eq1.2}
such that
$$\sum_{j=2}^{\infty}j(|a_{j}|+|b_{j}|)\leq 1-|b_{1}|,\quad (0\leq|b_{1}|<1),$$
and the subclass $HC$ of $HS$, where
$$\sum_{j=2}^{\infty}j^{2}(|a_{j}|+|b_{j}|)\leq 1-|b_{1}|,\quad (0\leq|b_{1}|<1).$$
The corresponding subclasses of $HS$ and $HC$ with $b_{1}=0$ are
denoted by $HS^{0}$ and $HC^{0}$, respectively. These two classes constitute a harmonic
counterpart of classes introduced by Goodman \cite{go}. They are useful in studying questions of so-called $\delta$-neighborhoods (Ruscheweyh \cite{ru}, see also \cite{qiwa}) and in constructing explicit $k$-quasiconformal extensions (Fait et al. \cite{fa}).

Our aim is to generalize the following result, due to Duren \cite{du}, to the mappings of the class
$HC_{p}^{0}$.

\begin{Thm}\label{ThmB} {\rm \bf (\cite[Theorem $1$, p. 50]{du})}
Each function $f\in C_{H}^{0}$ contains the full disk
$|w|<1/2$ in its range $f(\mathbb{U})$.\end{Thm}

A well-known coefficient conjecture of Clunie and Sheil-Small \cite{CPW4},
is that if $f=h+\overline{g}\in S_{H}^{0}$, then the Taylor coefficients
of the series of $h$ and $g$ satisfy the inequality
\be\label{eq00.1}|a_{j}|\leq \frac{1}{6}(2j+1)(j+1)\;\;\text{and} \;\;|b_{j}|\leq \frac{1}{6}(2j-1)(j-1)\ee
for all $j\geq1$. Although, this coefficient conjecture remains an open problem for the full class
$S_{H}^{0}$, this statement has been verified for certain subclasses, namely,
the class $T_{H}$ (see \cite[Section 6.6]{du}) of harmonic univalent typically real mappings, the class of
harmonic convex mappings in one direction, harmonic starlike mappings
in $S_{H}^{0}$ (see \cite[Section 6.7]{du}), and the class of harmonic
close-to-convex mappings (see \cite{wa}). Equality occurs in \eqref{eq00.1} for the harmonic Koebe mapping
\be\label{eq1.100}K(z)=\frac{z-\frac{1}{2}z^{2}+\frac{1}{6}z^{3}}{(1-z)^{3}}+
\overline{\frac{\frac{1}{2}z^{2}+\frac{1}{6}z^{3}}{(1-z)^{3}}},\ee
which is constructed by shearing the Koebe function $k(z) =z/(1-z)^{2}$
horizontally with the dilatation
$w(z) =z$. Note that $K$ maps the unit disk $\mathbb{U}$ onto the slit-plane $\mathbb{C}\setminus (-\infty,-1/6]$.

This paper is organized as follows. In Section 3, we generalize Theorem A to the class
$HC_{p}^{0}$ of polyharmonic mappings. The main result of this section is Theorem \ref{thm1.1}.
In Section 4, we obtain two convolution characterizations for polyharmonic
mappings to be starlike of order $\alpha$ and convex of order $\beta$, respectively. Our
results are Theorems \ref{thm2.1} and \ref{thm2.3},
where Theorem \ref{thm2.1} extends \cite[Theorems $2.6$]{ah},
and Theorem \ref{thm2.3} is a generalization of \cite[Theorem $2.8$]{ah}. In Section 5, we find the
radii of convexity and starlikeness for polyharmonic mappings, under certain coefficient conditions. The results in this section are Theorems $\ref{thm2.4}\sim\ref{thm2.7}$, which are the generalizations
of \cite[Theorems 3.1 and 3.3]{suna} and \cite[Theorem 1.11]{kal}, respectively.

\section{Preliminaries}\label{csw-sec2}

In this paper, we consider the polyharmonic mappings in
$\mathbb{U}$. We use $H_{p}$ to denote the set of all polyharmonic mappings $F$ in $\mathbb{U}$ with
a series expansion of the following form:

\be\label{eq1.1}
F(z)=\sum_{k=1}^{p}|z|^{2(k-1)}\big(h_{k}(z)+\overline{g_{k}(z)}\big)
=\sum_{k=1}^{p}|z|^{2(k-1)}\sum_{j=1}^{\infty}(a_{k,j}z^{j}+\overline{b_{k,j}}\overline{z^{j}}), \ee
with $a_{1,1}=1$, $|b_{1,1}|<1$. Let $H^{0}_{p}$
denote the subclass of $H_{p}$ for $b_{1,1}=0$
and $a_{k,1}=b_{k,1}=0$ for $k\in \{2,\cdots,p\}$.

\begin{defn}{\rm \bf (\cite{qiwa})}\label{defn1}
We say that a univalent polyharmonic mapping $F$ with $F(0) = 0$ is
starlike with respect to the origin if the curve $F(re^{i\theta})$
is starlike with respect to the
origin for each $r\in(0,1)$.
\end{defn}

\begin{prop}{\rm \bf (\cite{po})}\label{pro1}
If $F$ is univalent, $F(0) = 0$ and
$\frac{\partial}{\partial\theta}\big(\arg F(re^{i\theta})\big) >0$ for
$z = re^{i\theta}\not= 0$, then $F$ is starlike with respect to the origin.
\end{prop}

\begin{defn}{\rm \bf (\cite{qiwa})}\label{defn2}
A univalent polyharmonic mapping $F$ with $F(0) = 0$ and $\frac{\partial}{\partial \theta}F(re^{i\theta})$ $
\neq0$ whenever $r\in(0,1)$, is said to be convex
if the curve $F(re^{i\theta})$ is convex for each $r\in(0,1)$.
\end{defn}

\begin{prop}{\rm \bf (\cite{po})} \label{pro2}
 If $F$ is univalent,
$F(0) = 0$, $\frac{\partial}{\partial \theta}F(re^{i\theta})
\neq0$ whenever $r\in(0,1)$, and $\frac{\partial}{\partial \theta}
 \left[ \arg \left(\frac{\partial}{\partial \theta}F(re^{i\theta})\right) \right] >0$ for
$z = re^{i\theta}\not= 0$, then $F$ is convex.
\end{prop}

In \cite{qiwa}, J. Qiao and X. Wang introduced the subclass of $H_{p}^{0}$
denoted by $HS_{p}^{0}$ of polyharmonic mappings $F$ of the form \eqref{eq1.1} satisfying the condition
\be\label{eq1.5}
\sum_{k=1}^{p}\sum_{j=1}^{\infty}
\big(2(k-1)+j\big)\big(|a_{k,j}|+|b_{k,j}|\big)\leq2,\ee
and the subclass $HC_{p}^{0}$ of $HS_{p}^{0}$, where
\be\label{eq1.6}
\sum_{k=1}^{p}\sum_{j=1}^{\infty}
\big(2(k-1)+j^{2}\big)\big(|a_{k,j}|+|b_{k,j}|\big)\leq2.\ee
Their main result is the following:

\begin{Thm}{\rm \bf (\cite[Theorems $3.1$, $3.2$ and $3.3$]{qiwa})} \label{thmA}
Suppose $F\in HS_{p}^{0}$. Then $F$ is
univalent, sense preserving, starlike in $\mathbb{U}$.
In particularly, for each member of $HC_{p}^{0}$, $F$ maps
$\mathbb{U}$ onto a convex domain.
\end{Thm}

Obviously, if $p=1$, then the classes $HS^{0}_{p}$ and $HC^{0}_{p}$ reduce to $HS^{0}$ and $HC^{0}$, respectively.

\section{Coefficient estimates}\label{csw-sec3}
Now, we will generalize the Theorem A \cite{du} from the class $C_{H}^{0}$
to the class $HC_{p}^{0}$ of polyharmonic mappings.
\bthm\label{thm1.1} Let $F\in HC_{p}^{0}$ of the form \eqref{eq1.1}. Then the range $F(\mathbb{U})$
contains the full disk $|w|<1/2$.\ethm
\bpf Let $F\in HC_{p}^{0}$, and let $r\in(0,1)$. Write
$$F_{r}(z)=z+\sum_{j=2}^{\infty}\left(\sum_{k=1}^{p}
a_{k,j}r^{2(k-1)}\right)z^{j}
+\sum_{j=2}^{\infty}
\left(
\sum_{k=1}^{p}\overline{b_{k,j}}r^{2(k-1)}\right)\overline{z^{j}},\; z\in \mathbb{U}.
$$
Then $F_{r}$ is harmonic. By the hypothesis and \eqref{eq1.6}, $F\in HC_{p}^{0}$, which implies
$$\sum_{j=2}^{\infty}j^{2}\left|\sum_{k=1}^{p}r^{2(k-1)}a_{k,j} \right|
+\sum_{j=2}^{\infty}j^{2}\left|\sum_{k=1}^{p}r^{2(k-1)}b_{k,j} \right|\leq1,$$
i.e., $F_{r}\in C_{H}^{0}.$ As in the proof of Theorem A, we see that
the range $F_{r}(\mathbb{U})$ is convex. Thus, if
$w\not \in F_{r}(\mathbb{U})$, a suitable rotation gives
$$\mbox{Re} \left\{e^{i\theta}\big(F_{r}(z)-w\big)\right\}>0,$$
for all $z\in \mathbb{U}$. But if
$F_{r}(z)=\sum_{k=1}^{p}r^{2(k-1)}\big(h_{k}(z)+\overline{g_{k}(z)}\big)$,
it follows that $\rm{Re}\{\varphi(z)\}>0$ for
\begin{align*}
\varphi(z)=&\mbox{Re}\left\{e^{i\theta}\left( \sum_{k=1}^{p}r^{2(k-1)}h_{k}(z)-w\right)+e^{-i\theta}\sum_{k=1}^{p}r^{2(k-1)}g_{k}(z)\right\}\\
=&\mbox{Re}\{c_{0}+c_{1}z+\cdots\},\\
\end{align*}
where $c_{0}=-e^{i\theta} w$ and $c_{1}=e^{i\theta}$. Following the proof of Theorem A, we get
$$1=|e^{i\theta}|=|c_{1}|\leq2|c_{0}|=2|-e^{i\theta} w|=2|w|,$$
or $|w|\geq1/2$. This proves the result.
\epf

\begin{figure}
\centering
\includegraphics[width=2.8in]{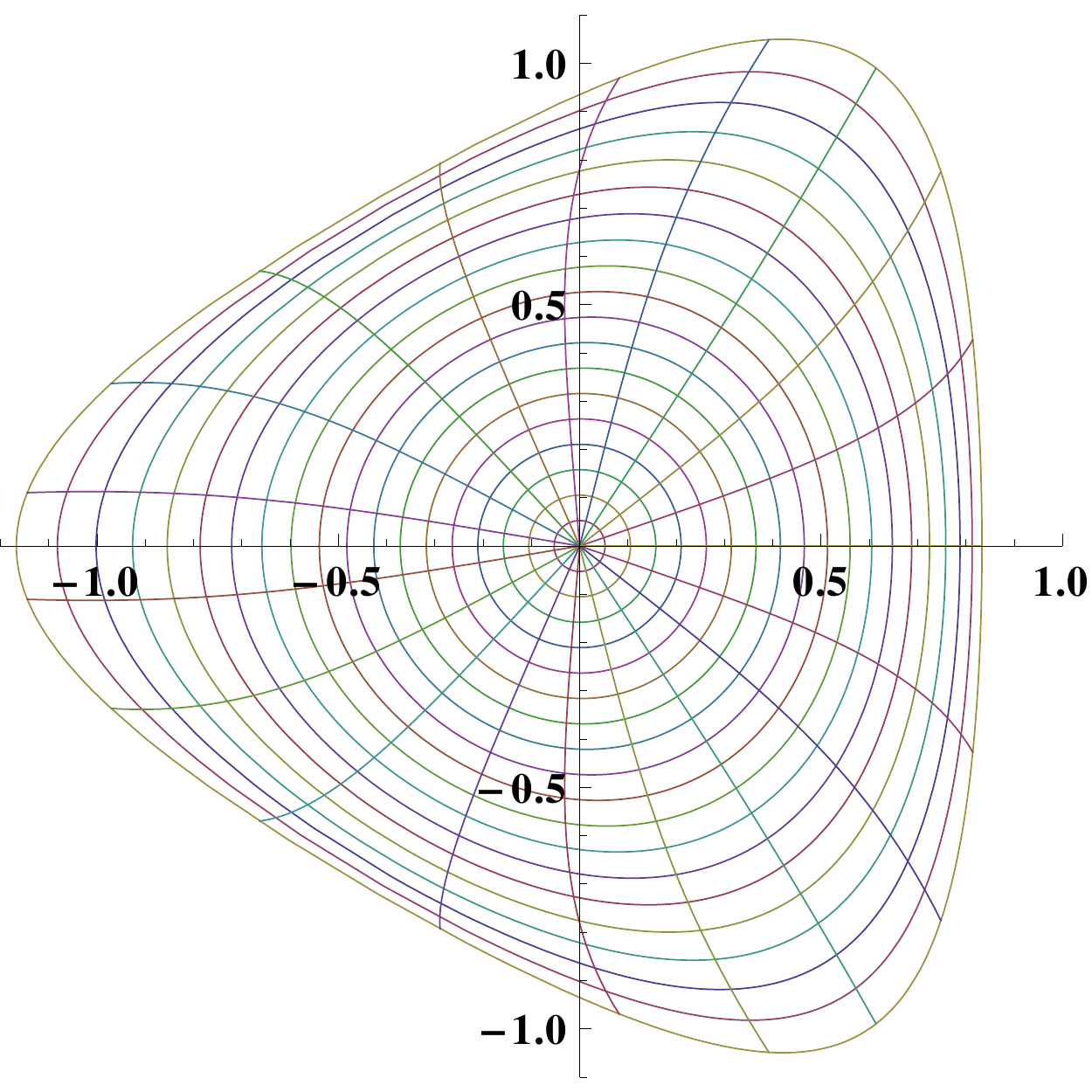}
\caption{The image of $\mathbb{U}$ under the mapping $F(z)=
z-\frac{1}{6}\overline{z^{2}}|z|^{2}.$ }\label{f1}
\end{figure}

\begin{examp} Let $F(z)=
z-\frac{1}{6}\overline{z^{2}}|z|^{2}\in HC_{2}^{0}$.
Then $F(\mathbb{U})$ contains the full disk $|w| < 1/2$. See Figure \ref{f1}.
\end{examp}

\section{Convolution characterization}\label{csw-sec4}
In this section, we obtain two convolution characterizations
concerning polyharmonic mappings which are starlike of order $\alpha$,
and convex of order $\beta$, respectively.

\begin{defn}{\rm \bf (\cite{luwa})}\label{defn3}
We say that a univalent polyharmonic mapping $F$ with $F(0) = 0$ is
starlike of order $\alpha\in[0,1)$ with respect to the origin if
\be\label{eq2.0}\frac{\partial}{\partial \theta}
\big(\arg F(re^{i\theta})\big)=\mbox{Re}
\left\{\frac{zF_{z}(z) -\overline{z}F_{\overline{z}}(z)}{F(z)}\right\} >\alpha\ee
for all $z=re^{i\theta}\neq0$.
\end{defn}

\begin{defn}{\rm \bf (\cite{luwa})}\label{defn4}
A univalent
polyharmonic mapping $F$ with $F(0) = 0$ and $\frac{\partial}{\partial \theta}F(re^{i\theta})$ $
\neq0$ whenever $r\in(0,1)$, is said to be convex of order $\beta\in[0,1)$
if \begin{align}\begin{split}\label{eq2.4}&\frac{\partial}{\partial \theta}
\left[\arg \left(\frac{\partial}{\partial \theta }F(re^{i\theta})\right)\right]\\
=&\mbox{Re}\left\{\frac{zF_{z}(z)
+z^{2}F_{z^{2}}(z)-2|z|^{2}F_{z\overline{z}}(z)+ \overline{z}F_{\overline{z}}(z)
+\overline{z}^{2}F_{\overline{z}^{2}}(z)}
{zF_{z}(z) -\overline{z}F_{\overline{z}}(z)}\right\} >\beta\end{split}\end{align}
for all $z=re^{i\theta}\neq0$.\end{defn}

\bthm\label{thm2.1} Let $F=\sum_{k=1}^{p}|z|^{2(k-1)}(h_{k}(z)+\overline{g_{k}(z)} )
\in H_{p}^{0}$ be univalent. Then $F$ is starlike of order $\alpha$ if and only if
\begin{align}\begin{split}\label{eq2.1}
&\sum_{k=1}^{p}|z|^{2(k-1)}
\left\{ h_{k}(z)\ast \left[\frac{z+((\alpha\xi+\alpha+\xi-1)/(2-\alpha-\alpha\xi))z^{2}}{(1-z)^{2}}\right]\right.\\
&\left.- \overline{g_{k}(z)}\ast \left[\frac{(2\xi+\alpha+\alpha\xi)/(2-\alpha-\alpha\xi) \overline{z}-((\alpha\xi+\alpha+\xi-1)/(2-\alpha-\alpha\xi))\overline{z}^{2}}
{(1-\overline{z})^{2}} \right] \right\}\neq 0,
\end{split}\end{align}
for all $z\neq 0$ in
$\mathbb{U}$ and all $\xi\in \IC$ with $|\xi|=1$.\ethm

\bpf Let $F\in H_{p}^{0}$ be univalent. Since
$$\frac{zF_{z}(z) -\overline{z}F_{\overline{z}}(z)}{F(z)}=1$$
at $z=0$, the condition \eqref{eq2.0} is equivalent to the condition
\beq\label{fri-1} \frac{zF_{z}(z) -\overline{z}F_{\overline{z}}(z)}{F(z)}-\alpha
\neq\frac{\xi-1}{\xi+1},\eeq for all $z\neq0$ in
$\mathbb{U}$ and all $\xi\in \mathbb{C}$ with $|\xi|=1$ and $\xi\neq-1$.
By the hypothesis that $F$ is univalent in $\mathbb{U}$, we get that $F(z)\not=0$
for $z\in \mathbb{U}\setminus\{0\}$.
Then, \eqref{fri-1} holds if and only if
\begin{align*}&(\xi+1)\left(\sum_{k=1}^{p}|z|^{2(k-1)}
\big(zh_{k}'(z)-\overline{z}\overline{g_{k}'(z})
-\alpha h_{k}(z)-\alpha\overline{g_{k}(z})\big)  \right)
\\
&\not=(\xi-1)\left(\sum_{k=1}^{p}|z|^{2(k-1)}
\big (h_{k}(z)+\overline{g_{k}(z})\big) \right)\end{align*}
for all $z\neq 0$ in
$\mathbb{U}$ and all $\xi\in \IC$ with $|\xi|=1$. Straightforward computations show that
\begin{align*}
&(\xi+1)\left(\sum_{k=1}^{p}|z|^{2(k-1)}
\big(zh_{k}'(z)-\overline{z}\overline{g_{k}'(z})
-\alpha h_{k}(z)-\alpha\overline{g_{k}(z})\big)  \right)\\
&-(\xi-1)\left(\sum_{k=1}^{p}|z|^{2(k-1)}
\big (h_{k}(z)+\overline{g_{k}(z})\big) \right)\\
=&\sum_{k=1}^{p}|z|^{2(k-1)}\left\{ h_{k}(z)\ast
 \left[\frac{(\xi+1)z}{(1-z)^{2}}- \frac{(\alpha\xi+\alpha+\xi-1)z}{1-z}  \right]\right.\\
 &\left.- \overline{g_{k}(z)}\ast \left[\frac{(\xi+1)\overline{z}}{(1-\overline{z})^{2}}
+\frac{(\alpha\xi+\alpha+\xi-1)\overline{z}}{1-\overline{z}}  \right] \right\} \\
=& \sum_{k=1}^{p}|z|^{2(k-1)}\left\{ h_{k}(z)\ast
\left[\frac{(2-\alpha-\alpha\xi)z+(\alpha\xi+\alpha+\xi-1)z^{2}}{(1-z)^{2}}\right]\right.\\
&\left.- \overline{g_{k}(z)}\ast \left[\frac{(2\xi+\alpha+\alpha\xi) \overline{z}-(\alpha\xi+\alpha+\xi-1)\overline{z}^{2}}{(1-\overline{z})^{2}} \right] \right\} ,\\
\end{align*}
from which we see that \eqref{fri-1} is true if and only if so is
\eqref{eq2.1}. The proof is complete. \epf

\br\label{cor2.2}
The above result gives a sufficient condition for mappings in $H_{p}^{0}$ to be starlike in terms of their coefficients. Let $F\in H_{p}^{0}$ be of the form \eqref{eq1.1}. If
$$\sum_{k=1}^{p}\sum_{j=2}^{\infty}\frac{2(k-1)+j-\alpha}{1-\alpha}|a_{k,j}|
+\sum_{k=1}^{p}\sum_{j=2}^{\infty}\frac{2(k-1)+j+\alpha}{1-\alpha}|b_{k,j}|\leq1,
$$ then $F$ is sense-preserving, univalent and starlike of order $\alpha$.
The result follows from Theorem \ref{thm2.1} and Lemma B by a straightforward calculation. In fact, this case is already covered by Theorem B.\er

\bthm\label{thm2.3} Let $F=\sum_{k=1}^{p}|z|^{2(k-1)}
\big(h_{k}(z)+\overline{g_{k}(z)}\big)  \in H_{p}^{0}$
be univalent such that $\frac{\partial}{\partial\theta}F(re^{i \theta})\neq0$
for all $r\in (0,1)$. Then $F$ is convex of order $\beta$ if and only if
\begin{align}\begin{split}\label{eq4.3}
&\sum_{k=1}^{p}|z|^{2(k-1)} \left\{h_{k}(z)
\ast\left[\frac{(2-\beta\xi-\beta)z+(2\xi+\beta\xi+\beta) z^{2}}{(1-z)^{3}}\right]\right. \\
&\left.+ \overline{g_{k}(z)}  \ast \left[\frac{(2\xi+\beta\xi+\beta) \overline{z}+(2-\beta\xi-\beta)
\overline{z}^{2}}{(1-\overline{z})^{3}} \right]\right\}\neq 0,
\end{split}\end{align}
for all $z\neq 0$ in
$\mathbb{U}$ and all $\xi\in \IC$ with $|\xi|=1$.\ethm

\bpf Let $F\in H_{p}^{0}$ be univalent. Since $$\frac{zF_{z}(z)+z^{2}F_{z^{2}}(z)-2|z|^{2}F_{z\overline{z}}(z)+
\overline{z}F_{\overline{z}}(z)+\overline{z}^{2}F_{\overline{z}^{2}}(z) }
{zF_{z}(z) -\overline{z}F_{\overline{z}}(z)} =1$$ at $z=0$, the required condition \eqref{eq2.4} is equivalent to
\beq\label{fri-2} \frac{zF_{z}(z)+z^{2}F_{z^{2}}(z)-2|z|^{2}F_{z\overline{z}}(z)+ \overline{z}F_{\overline{z}}(z)+\overline{z}^{2}F_{\overline{z}^{2}}(z)}{zF_{z}(z) -\overline{z}F_{\overline{z}}(z)}-\beta\neq\frac{\xi-1}{\xi+1},\eeq for all $z\neq0$ in
$\mathbb{U}$ and all $\xi\in \mathbb{C}$ with $|\xi|=1$ and $\xi\neq-1$.
Note that $\frac{\partial}{\partial \theta}F(re^{i\theta})
\not=0$ for all $r\in(0,1)$.
Then, \eqref{fri-2} holds if and only if
\begin{align*}&(\xi+1)\sum_{k=1}^{p}|z|^{2(k-1)}\big((1-\beta)zh'_{k}(z)+z^{2}h''_{k}(z)+ (1+\beta)\overline{z}\overline{g'_{k}(z)}+\overline{z}^{2}\overline{g''_{k}(z)}\big)\\
&-(\xi-1)\sum_{k=1}^{p}|z|^{2(k-1)}\big(zh'_{k}(z) -\overline{z}\overline{g'_{k}(z)} \big)\not=0
\end{align*}
for all $z\neq 0$ in
$\mathbb{U}$ and all $\xi\in \IC$ with $|\xi|=1$. Straightforward computations show that

\begin{align*}
&(\xi+1)\sum_{k=1}^{p}|z|^{2(k-1)}\big((1-\beta)zh'_{k}(z)+z^{2}h''_{k}(z)+ (1+\beta)\overline{z}\overline{g'_{k}(z)}+\overline{z}^{2}\overline{g''_{k}(z)}\big)\\
&-(\xi-1)\sum_{k=1}^{p}|z|^{2(k-1)}\big(zh'_{k}(z) -\overline{z}\overline{g'_{k}(z)} \big)\\
=&\sum_{k=1}^{p}|z|^{2(k-1)}\left\{ h_{k}(z)\ast \left[\frac{z(2-\beta\xi-\beta)}{(1-z)^{2}}
+\frac{2z^{2}(\xi+1)}{(1-z)^{3}}  \right]\right.\\
&\left.+ \overline{g_{k}(z)}\ast \left[\frac{\overline{z}(2\xi+\beta\xi+\beta)}{(1-\overline{z})^{2}}
+\frac{2\overline{z}^{2}(\xi+1)}{(1-\overline{z})^{3}} \right] \right\} \\
=& \sum_{k=1}^{p}|z|^{2(k-1)} \left\{h_{k}(z)\ast
\left[\frac{(2-\beta\xi-\beta)z+(2\xi+\beta\xi+\beta) z^{2}}{(1-z)^{3}}\right]\right.\\
 &\left.+ \overline{g_{k}(z)}  \ast \left[\frac{(2\xi+\beta\xi+\beta) \overline{z}+
 (2-\beta\xi-\beta)\overline{z}^{2}}{(1-\overline{z})^{3}} \right]\right\} ,\\
\end{align*}
from which we see that \eqref{fri-2} is true if and only if \eqref{eq4.3} is. The proof is complete. \epf

\br\label{cor2.3}
By a straightforward calculation, we obtain from Theorem \ref{thm2.3} and Lemma B a sufficient coefficient bound for polyharmonic mappings which are convex of order $\beta$. Let $F\in H_{p}^{0}$ be of the form \eqref{eq1.1}. If
$$\sum_{k=1}^{p}\sum_{j=2}^{\infty}\frac{2(k-1)+j^{2}-\beta}{1-\beta}|a_{k,j}|
+\sum_{k=1}^{p}\sum_{j=1}^{\infty}\frac{2(k-1)+j^{2}+\beta}{1-\beta}|b_{k,j}|\leq1,
$$
then $F$ is convex of order $\beta$. In fact, this case is already covered by Theorem B.
\er

\section{Radii for starlikeness and convexity}\label{csw-sec5}
In this section, we will first generalize the results \cite[Theorems $3.1$ and $3.3$]{suna} to the polyharmonic mappings. The following identities, where $r\in(0,1)$, are used in the proofs of our results:
\begin{align}\begin{split}\label{eq00002.4}
&\sum_{j=1}^{\infty}r^{j-1}=\frac{1}{1-r},\;\;\sum_{j=1}^{\infty}jr^{j-1}=\frac{1}{(1-r)^{2}},\;
\;\sum_{j=1}^{\infty}j^{2}r^{j-1}=\frac{1+r}{(1-r)^{3}},\\
&\sum_{j=1}^{\infty}j^{3}r^{j-1}=\frac{1+4r+r^{2}}{(1-r)^{4}}\;\;\text{and}\;
\;\sum_{j=1}^{\infty}j^{4}r^{j-1}=\frac{(1+r)(1+10r+r^{2})}{(1-r)^{5}}.\\
\end{split}\end{align}

\bthm\label{thm2.4} Let $F\in H_{p}^{0}$ of the form \eqref{eq1.1} and the
coefficients of the series satisfy the conditions
\begin{center}
$|a_{k,j}|\leq\frac{1}{6}(2j+1)(j+1)$ and $|b_{k,j}|\leq\frac{1}{6}(2j-1)(j-1).$
\end{center}
 Then $F$ is univalent and starlike of order $\alpha$ in $|z|<r_{0}(\alpha)$,
 where $r_{0}(\alpha)$ is the smallest positive root of the equation
\be\label{eq2.3} 6(1-\alpha)(1-r)^{4}-
\sum_{k=1}^{p}r^{2(k-1)}\big(3(r+1)^{2}-3\alpha(1-r)^{2}+2(k-1)(r^{2}+3)(1-r)  \big)=0,\ee
in the interval $(0,1)$. The result is sharp.
\ethm
\bpf
Let $F_{r}(z):=r^{-1}F(rz)$, where $F\in H_{p}^{0}$ is of the form \eqref{eq1.1}, and fix $r\in(0,1)$. Then
$$F_{r}(z)=\sum_{k=1}^{p}|z|^{2(k-1)}\sum_{j=1}^{\infty}
\left(a_{k,j}r^{2k+j-3}z^{j}+\overline{b_{k,j}}r^{2k+j-3}
\overline{z^{j}}\right),\; z\in \mathbb{U}.$$
By the hypotheses, $|a_{k,j}|\leq\frac{1}{6}(2j+1)(j+1)$ and $|b_{k,j}|
\leq\frac{1}{6}(2j-1)(j-1)$. By using these coefficient estimates and \eqref{eq00002.4}, we obtain
\begin{align*}
S_0:=&\sum_{k=1}^{p}\sum_{j=1}^{\infty}\left(\frac{2(k-1)+j-\alpha}{1-\alpha}|a_{k,j}|r^{2k+j-3}
+\frac{2(k-1)+j+\alpha}{1-\alpha}|b_{k,j}|r^{2k+j-3}\right)\\
\leq&\frac{1}{6}\sum_{k=1}^{p}\sum_{j=1}^{\infty}
\left(\frac{2(k-1)+j-\alpha}{1-\alpha}(2j+1)(j+1)\right.\\
&\left.+\frac{2(k-1)+j+\alpha}{1-\alpha}(2j-1)(j-1)\right)r^{2k+j-3}\\
=&\sum_{k=1}^{p}\sum_{j=1}^{\infty}\frac{1}{3(1-\alpha)}
\big(2j^{3}+4j^{2}(k-1)+(1-3\alpha)j+2(k-1)\big)r^{2k+j-3}\\
=&\frac{1}{3(1-\alpha)}\sum_{k=1}^{p}\left(2r^{2k-3}
\sum_{j=1}^{\infty}j^{3}r^{j}+4(k-1)r^{2k-3}\sum_{j=1}^{\infty}j^{2}r^{j}\right.\\
&\left.+(1-3\alpha)r^{2k-3}\sum_{j=1}^{\infty}jr^{j}
+2(k-1)r^{2k-3}\sum_{j=1}^{\infty}r^{j}  \right)\\
=&\frac{1}{3(1-\alpha)}\sum_{k=1}^{p}r^{2k-2}\left(\frac{3(r+1)^{2}}{(1-r)^{4}}
-\frac{3\alpha}{(1-r)^{2}}+\frac{2(k-1)(r^{2}+3)}{(1-r)^{3}} \right).\\
\end{align*}
According to Remark 1, it suffices to show that $S_0\leq2$.
By the last inequality, $S_0\leq2$ if $r$ satisfies the inequality
\begin{align*}s_0(r)=&6(1-\alpha)(1-r)^{4}-\sum_{k=1}^{p}r^{2(k-1)}\big(3(r+1)^{2}
-3\alpha(1-r)^{2}\\
&+2(k-1)(r^{2}+3)(1-r) \big)\geq0.
\end{align*}
Since $s_0(0)=3-3\alpha>0$ and $s_0(1)<0$, then exists a smallest positive
root $r_{0}(\alpha)$ of the equation $s_0(r)=0$ in the interval (0,1). In
particular, $F$ is sense-preserving, univalent and starlike of order $\alpha$ in $|z|<r_{0}(\alpha)$.

As in \cite[Theorem $3.1$]{suna}, the mapping
\begin{equation}
\label{f0def}
f_{0}(z)=2z-\frac{z-\frac{1}{2}z^{2}+\frac{1}{6}z^{3}}{(1-z)^{3}}
 +\overline{\frac{\frac{1}{2}z^{2}+\frac{1}{6}z^{3}}{(1-z)^{3}}}
\end{equation}
shows that the bound given by $r_{0}(\alpha)$ is the best possible.
\epf

\begin{figure}
\centering
\includegraphics[width=2.8in]{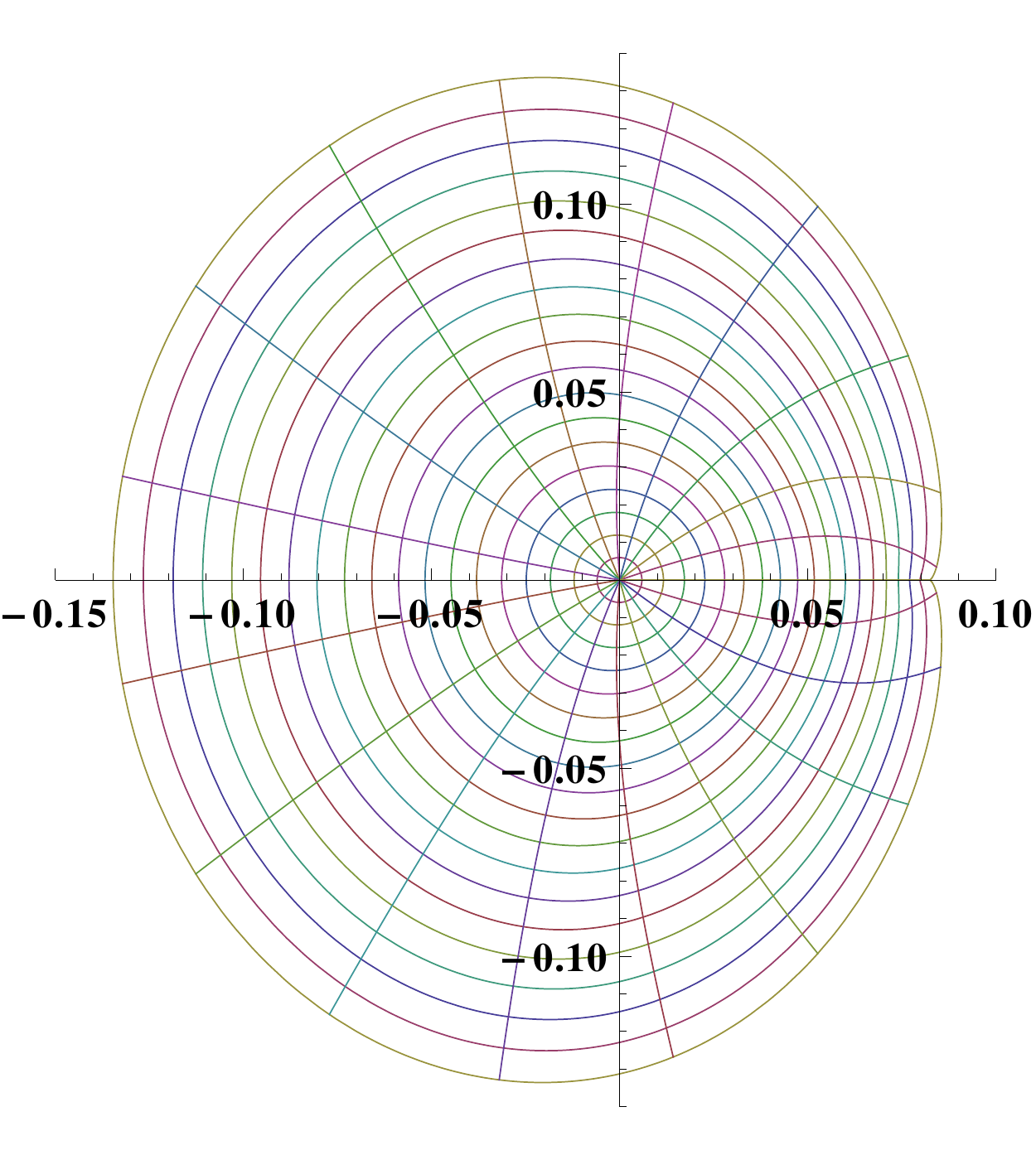}
\caption{The image of $\mathbb{U}(r_0(\alpha))$, where $\alpha=0$ and $r_0\approx 0.11290$ is given by \eqref{eq2.3}, under the mapping $f_0$ of \eqref{f0def}.}\label{f0fig}
\end{figure}

\begin{figure}
\centering
\includegraphics[width=2.8in]{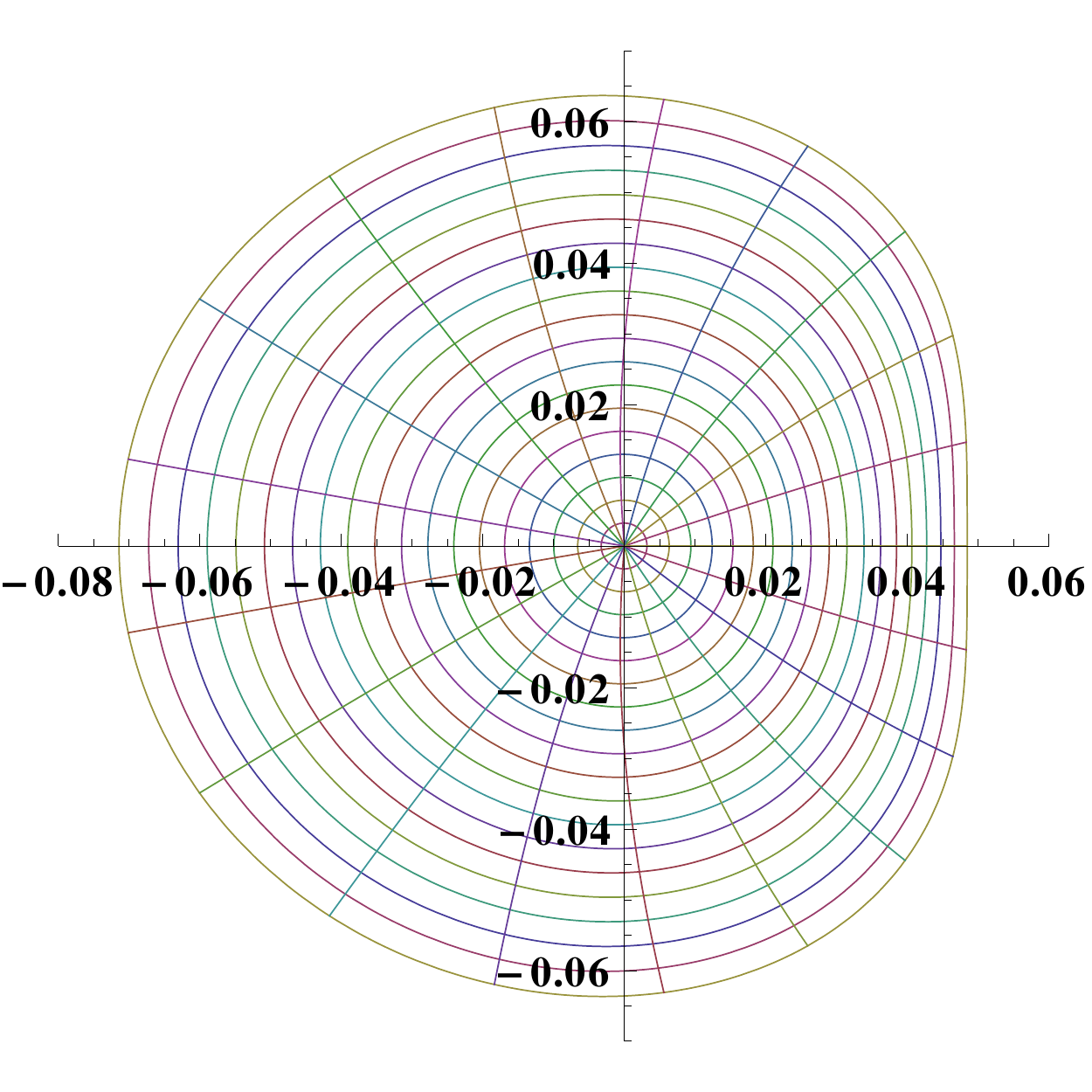}
\caption{The image of $\mathbb{U}(r_1(\beta))$, where $\beta=0$ and $r_1\approx 0.06143$ is given by \eqref{eq3.5}, under the mapping $f_1$ of \eqref{f1def}.}\label{f1fig}
\end{figure}

\bthm\label{thm2.5} Under the hypothesis of
Theorem \ref{thm2.4}, $F\in H_{p}^{0}$ is univalent and convex of
order $\beta$ in the disk $|z|<r_{1}(\beta)$, where
$r_{1}(\beta)$ is the smallest positive root of the equation
\begin{align}\begin{split}\label{eq3.5}
0=&6(1-\beta)(1-r)^{5}-\sum_{k=1}^{p}r^{2(k-1)}
\big((8k-6-6\beta)(1+r)(1-r)^{2}\\
&+4(k-1)(1-r)^{4}+4(1+r)(1+10r+r^{2})-6(2k-1-\beta)(1-r)^{5}\big)
\end{split}\end{align}
in the interval $(0,1)$. The result is sharp.
\ethm
\bpf The proof of this result is similar to Theorem \ref{thm2.4},
where Remark \ref{cor2.3} is used instead of Remark \ref{cor2.2}, and we omit it.
The bound $r_{1}(\beta)$ given by (\ref{eq3.5}) is again sharp, which can be seen by considering the mapping
\begin{equation}\label{f1def}
f_{1}(z)=2z-\frac{z-\frac{1}{2}z^{2}+\frac{1}{6}z^{3}}{(1-z)^{3}}
 -\overline{\frac{\frac{1}{2}z^{2}+\frac{1}{6}z^{3}}{(1-z)^{3}}},
\end{equation}
see Figure \ref{f1fig}. Then, by \cite[Theorem 3.3]{suna}, we see that the bound given by $r_{1}(\beta)$ is the best possible.
\epf

In \cite[Theorem $1.11$]{kal}, D. Kalaj et al. gave the radius of close-to-convexity of harmonic mappings under certain coefficients conditions. Now, we will study the radius of starlikeness and convexity of mappings in $H_{p}^{0}$ under the same coefficients condition. Our result is the following:

\bthm\label{thm2.6} Let $F\in H_{p}^{0}$ of the form \eqref{eq1.1} and the
coefficients of the series satisfy the conditions
\begin{center}
$|a_{k,j}|+|b_{k,j}|\leq C$ for all $j\geq2$.
\end{center}
 Then $F$ is univalent and starlike of order $\alpha$ in $|z|<r_{2}(\alpha)$,
 where $r_{2}(\alpha)$ is the smallest positive root of the equation
\be\label{eq2.5} (1-\alpha)(1-r)^{2}-
\sum_{k=1}^{p}Cr^{2(k-1)}\big((2k-2+\alpha)(1-r)+1-(2k+\alpha-1)(1-r)^{2}\big)=0\ee
in the interval (0,1). The result is sharp.
\ethm
\bpf Let $F_{r}(z):=r^{-1}F(rz)$, where $F\in H_{p}^{0}$ is of the form \eqref{eq1.1}, and fix $r\in(0,1)$. Then
$$F_{r}(z)=\sum_{k=1}^{p}|z|^{2(k-1)}\sum_{j=1}^{\infty}
\left(a_{k,j}r^{2k+j-3}z^{j}+\overline{b_{k,j}}r^{2k+j-3}
\overline{z^{j}}\right),\; z\in \mathbb{U}.$$
As in the proof of Theorem \ref{thm2.4}, under the hypothesis that $|a_{k,j}|+|b_{k,j}|
\leq C$ and \eqref{eq00002.4}, we get
\begin{align*}
S_1:=&\sum_{k=1}^{p}\sum_{j=2}^{\infty}\left(\frac{2(k-1)+j-\alpha}{1-\alpha}|a_{k,j}|
+\frac{2(k-1)+j+\alpha}{1-\alpha}|b_{k,j}|\right)r^{2k+j-3}\\
\leq&\sum_{k=1}^{p}\sum_{j=2}^{\infty}
\frac{2(k-1)+j+\alpha}{1-\alpha}Cr^{2k+j-3}\\
=&\sum_{k=1}^{p}\left(\frac{(2k-2+\alpha)(1-r)+1}{(1-r)^{2}} -2k+1-\alpha \right)\frac{Cr^{2k-2}}{1-\alpha}.
\end{align*}
According to Remark 1, it suffices to show that $S_1\leq1$.
By the last inequality, $S_1\leq1$ if $r$ satisfies the following inequality:
$$s_1(r)=(1-\alpha)(1-r)^{2}-
\sum_{k=1}^{p}Cr^{2(k-1)}\big((2k-2+\alpha)(1-r)+1-(2k+\alpha-1)(1-r)^{2}\big)\geq0.$$
Since $s_1(0)=1-\alpha>0$ and $s_1(1)<0$, then exists a smallest positive
root $r_{2}(\alpha)$ of the equation $s_1(r)=0$ in the interval $(0,1)$. In
particular, $F$ is sense-preserving, univalent and starlike of order $\alpha$ in $|z|<r_{2}(\alpha)$.

\begin{figure}
\centering
\includegraphics[width=2.8in]{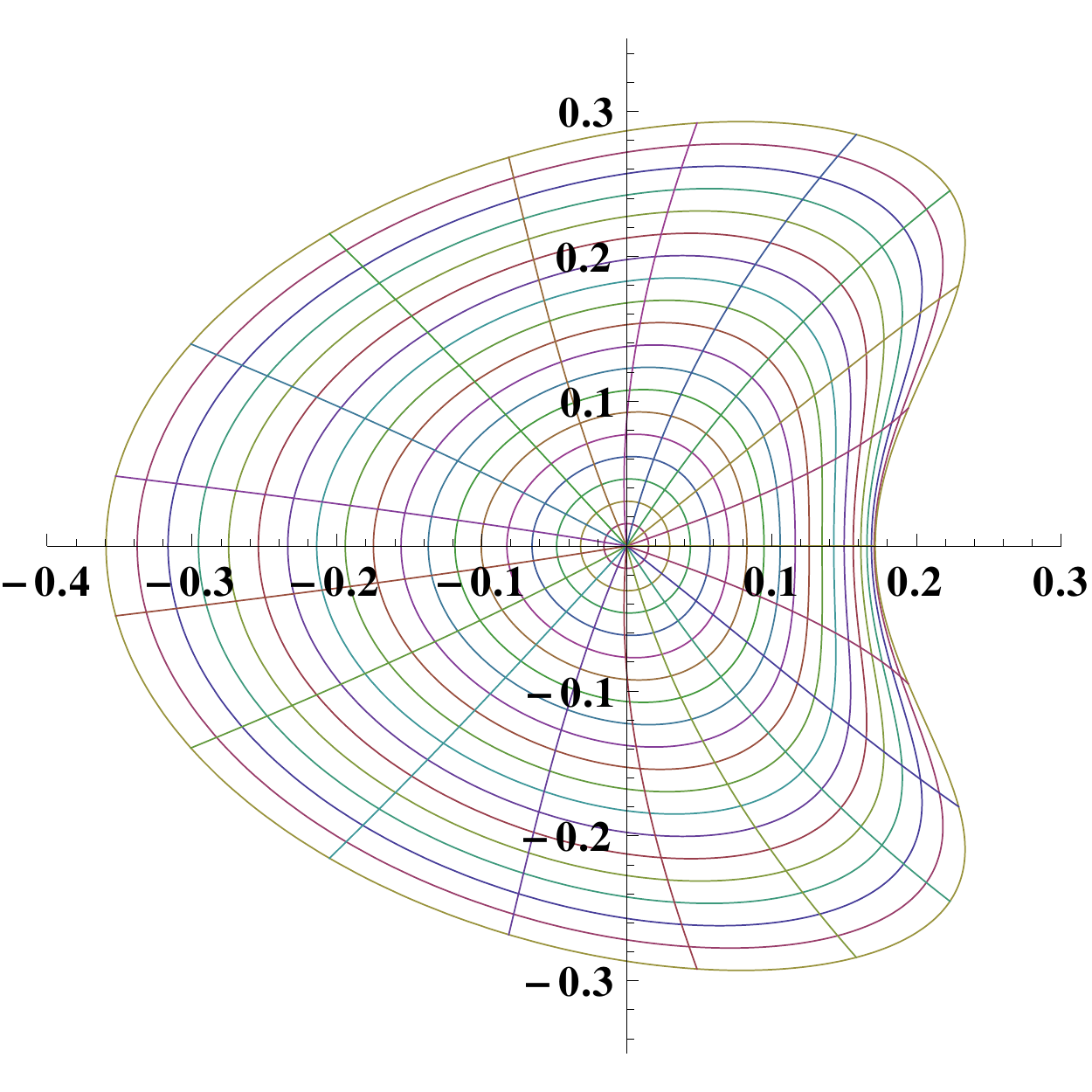}
\caption{The image of $\mathbb{U}(r_2(\alpha))$, where $\alpha=0$, $C=1$ and $r_2\approx 0.29289$ is given by \eqref{eq2.5}, under the mapping $f_2$ of \eqref{f2def}.}\label{f2fig}
\end{figure}

To prove the sharpness part of the statement, one may consider the mapping
\begin{equation}
\label{f2def}
f_{2}(z)=z-\frac{Cz^{2}}{2(1-z)}-\overline{ \frac{Cz^{2}}{2(1-z)}  },
\end{equation}
see Figure \ref{f2fig}.
Then by \cite[Theorem 1.11]{kal}, we see that the bound given by $r_{2}(\alpha)$ is the best possible. The proof of the theorem is complete.
\epf

\begin{figure}
\centering
\includegraphics[width=2.8in]{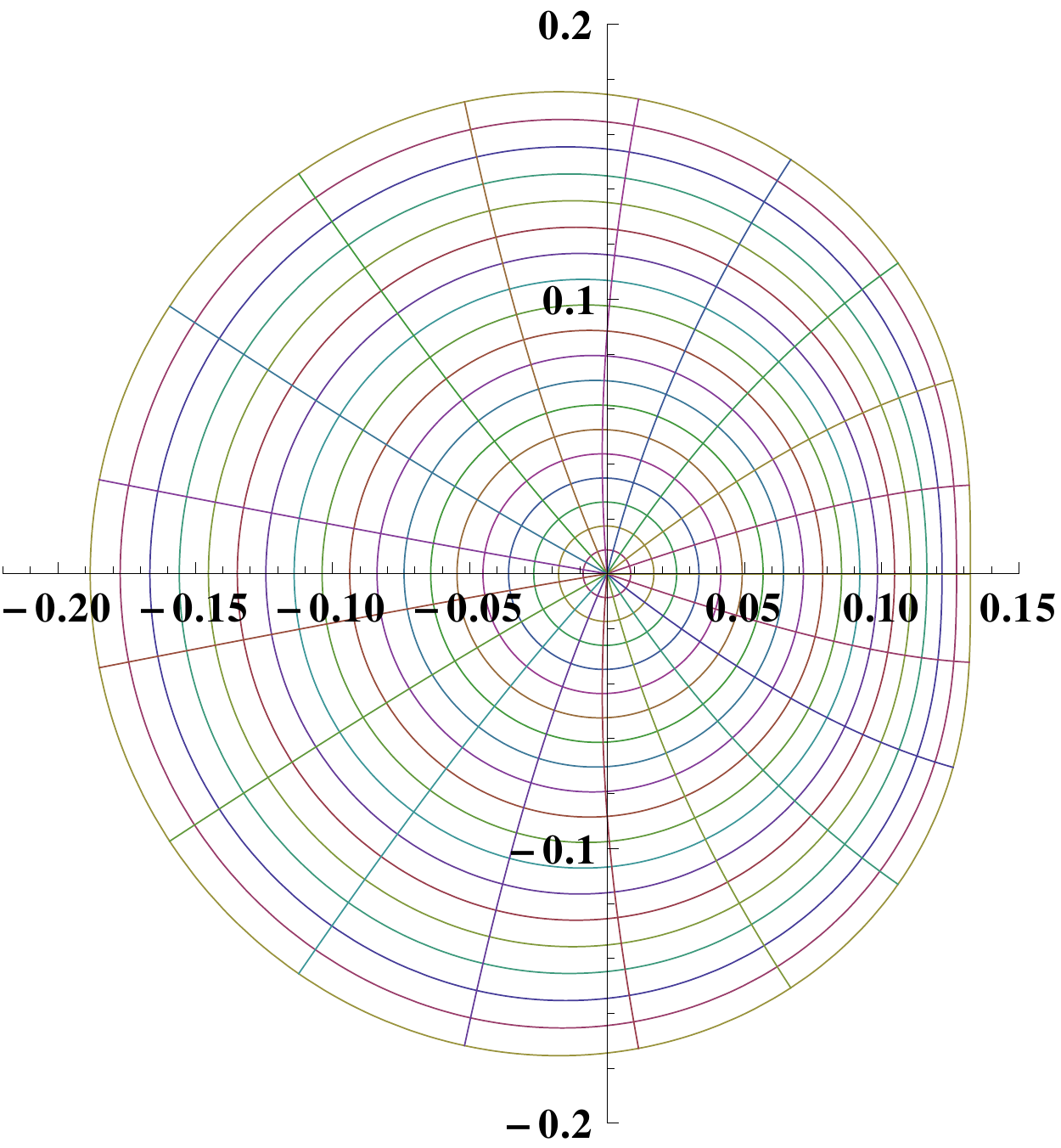}
\caption{The image of $\mathbb{U}(r_3(\beta))$, where $\beta=0$, $C=1$ and $r_3\approx 0.16488$ is given by \eqref{eq3.8}, under the mapping $f_3$ of \eqref{f3def}.}\label{f3fig}
\end{figure}

\bthm\label{thm2.7} Under the hypothesis of
Theorem \ref{thm2.6}, $F \in H_{p}^{0}$ is univalent and convex of order $\beta$ in the disk $|z|<r_{3}(\beta)$, where
$r_{3}(\beta)$ is the smallest positive real root of the equation
\begin{align}\begin{split}\label{eq3.8}
0=&(1-\beta)(1-r)^{3}-
\sum_{k=1}^{p}Cr^{2k-2}\big((2k-2)(1-r)^{2}+1+r+\beta-\beta r)\\
&-(2k+\beta-1)(1-r)^{3}\big)
\end{split}\end{align}
in the interval $(0,1)$. The result is sharp.
\ethm
\bpf The proof of this result is similar to that of Theorem \ref{thm2.6}, where Remark \ref{cor2.3} is used instead of Remark \ref{cor2.2}, and we omit it.
The bound $r_{3}(\beta)$ given by (\ref{eq3.8}) is sharp by considering the mapping
\begin{equation}
\label{f3def}
f_{3}(z)=z-\overline{\frac{Cz^{2}}{1-z}},
\end{equation}
see Figure \ref{f3fig}. Note that the root of the equation \eqref{eq3.8} in $(0,1)$ is decreasing as a function of $\beta\in[0,1)$. As $f_{3}$ has real coefficients, we obtain
$$\frac{\partial}{\partial \theta}
\left[\arg \left(\frac{\partial}{\partial \theta }f_{3}(re^{i\theta})\right)\right] \bigg|_{\theta=0,r=r_{3}(\beta)}=\frac{(1+C)-(C+Cr)/(1-r)^{3}}{1-C+C/(1-r)^{2}} \bigg|_{r=r_{3}(\beta)}=\beta. $$ Therefore, the mapping
$f_{3} $ will not be convex of order $\beta$ in the disk $|z|<r$, where $r>r_{3}(\beta)$.
\epf

\begin{thebibliography}{HD}




\normalsize
\baselineskip=17pt

%

\bibitem{abab} {\sc Z. Abdulhadi} and {\sc Y. Abu Muhanna}, Landau's theorem for biharmonic mappings.
{\it J. Math. Anal. Appl.} {\bf 338} (2008), 705--709.

\bibitem{A} {\sc Z. Abdulhadi}, {\sc Y. Abu Muhanna} and {\sc S. Khuri}, On univalent solutions of the
biharmonic equation. {\it J. Inequal. Appl.} {\bf 5} (2005), 469--478.

\bibitem{ababkh} {\sc Z. Abdulhadi}, {\sc Y. Abu Muhanna} and {\sc S. Khuri},
On some properties of solutions of the biharmonic equation.
{\it Appl. Math. Comput.} {\bf 117} (2006), 346--351.

\bibitem{ah} {\sc O. P. Ahuja} and {\sc J. M. Jahangiri},
Convolutions for special classes of harmonic univalent functions.
{\it Appl. Math. Lett.} {\bf 16} (2003), 905--909.

\bibitem{av} {\sc Y. Avci} and {\sc E. Z{\l}otkiewicz},
On harmonic univalent mappings.
{\it Ann. Univ. Mariae Curie Sk{\l}odowska (Sect A)} {\bf 44} (1990), 1--7.

\bibitem{CPW0} {\sc Sh. Chen}, {\sc S. Ponnusamy} and {\sc X. Wang},
Landau's theorem for certain biharmonic mappings.
{\it Appl. Math. Comput.} {\bf 208} (2009), 427--433.

\bibitem{CPW4} {\sc Sh. Chen}, {\sc S. Ponnusamy} and {\sc X. Wang},
Compositions of harmonic mappings and biharmonic
mappings. {\it Bull. Belg. Math. Soc. Simon
Stevin.} {\bf 17} (2010), 693--704.

\bibitem{CRW} {\sc J. Chen}, {\sc A. Rasila} and {\sc X. Wang}, On polyharmonic univalent mappings. {\it arXiv:1302.2018}

\bibitem{CW} {\sc J. Chen} and {\sc X. Wang}, On certain classes of biharmonic mappings
defined by convolution. {\it Abstr. Appl. Anal.}.
{\bf 2012}, Article ID 379130, 10 pages. doi:10.1155/2012/379130

\bibitem{CPW2} {\sc Sh.~Chen}, {\sc S.~Ponnusamy} and {\sc X.~Wang},
Bloch constant and Landau's theorems for planar $p$-harmonic
mappings. {\it J. Math. Anal. Appl.} {\bf 373} (2011), 102--110.

\bibitem{cl} {\sc J. G. Clunie} and {\sc T. Sheil-Small},
Harmonic univalent functions. {\it Ann. Acad. Sci. Fenn. Ser. A. I.} {\bf 9} (1984), 3--25.

\bibitem{du} {\sc P. Duren}, {\it Harmonic mappings in the plane}.
Cambridge University Press, Cambridge, 2004.

\bibitem{fa} {\sc M. Fait, J. Krzy\.{z}} and {\sc J. Zygmunt},
Explicit quasiconformal extensions for some classes of univalent functions.
{\it Comment. Math. Helv.} {\bf 51} (1976), 279--285.

\bibitem{go} {\sc A. W. Goodman}, Univalent functions and nonanalytic curves,
{\it Proc. Amer. Math. Soc.} {\bf 8} (1957), 588--601.

\bibitem{ha} {\sc J. Happel} and {\sc H. Brenner},
{\it Low Reynolds Number Hydrodynamics with Special Applications to Particulate Media}.
Prentice-Hall, Englewood Cliffs, NJ, USA,
1965.

\bibitem{kal} {\sc D. Kalaj}, {\sc S. Ponnusamy} and {\sc M. Vuorinen},
Radius of close-to-convexity of harmonic funcions.
{\it arXiv:1107.0610}.

\bibitem{kh} {\sc S. A. Khuri},
Biorthogonal series solution of Stokes flow problems in sectorial
regions. {\it SIAM J. Appl. Math.} {\bf 56} (1996), 19--39.

\bibitem{la} {\sc W. E. Langlois}, {\it Slow Viscous Flow}. Macmillan, New York, NY, USA, 1964.

\bibitem{luwa} {\sc Q. Luo} and {\sc X. Wang}, The starlikeness, convexity, covering theorem and extreme
points of $p$-harmonic mappings.
{\it Bull. Iranian Math. Soc.}, in press. 

\bibitem{suna} {\sc S. Nagpal} and {\sc V. Ravichandran},
Fully starlike and convex harmonic mappings of order $\alpha$.
{\it arXiv:1207.3946}.

\bibitem{qiwa} {\sc J. Qiao} and {\sc X. Wang},
On $p$-harmonic univalent mappings (in Chinese).
{\it Acta Math. Sci.} {\bf 32A} (2012), 588--600.

\bibitem{po} {\sc C. Pommerenke}, {\it Univalent functions}.
Vandenhoeck and Ruprecht, G\"{o}ttin-gen, 1975.

\bibitem{wa} {\sc X. Wang} and {\sc X. Liang},
Precise coefficient estimates for Close-to-convex harmonic univalent mappings.
{\it J. Math. Anal. Appl.} {\bf 263} (2001), 501--509.

\bibitem{ru} {\sc S. Ruscheweyh}, Neighborhoods of univalent functions.
{\it Proc. Amer. Math. Soc.} {\bf 18} (1981), 521--528.
\end{thebibliography}
\end{document}